\documentclass[11pt]{amsart}
\usepackage{times,epsf,amssymb,amsmath}

\begin{document}

\newtheorem{thm}{Theorem}[section]
\newtheorem{lem}[thm]{Lemma}
\newtheorem{cor}[thm]{Corollary}

\theoremstyle{definition}
\newtheorem{defn}{Definition}[section]

\theoremstyle{remark}
\newtheorem{rmk}{Remark}[section]

\def\square{\hfill${\vcenter{\vbox{\hrule height.4pt \hbox{\vrule
width.4pt height7pt \kern7pt \vrule width.4pt} \hrule height.4pt}}}$}

\def\T{\mathcal T}

\newenvironment{pf}{{\it Proof:}\quad}{\square \vskip 12pt}

\title{Properly Embedded Least Area Planes in Gromov Hyperbolic $3$-Spaces}
\author{Baris Coskunuzer}
\address{Department of Mathematics \\ Yale University \\ New Haven, CT 06520}
\email{baris.coskunuzer@yale.edu}
\thanks{The author is supported by NSF Grant DMS-0603532}

\maketitle


\newcommand{\cirD}{\overset{\circ}{D}}
\newcommand{\Si}{S^2_{\infty}(X)}
\newcommand{\SH}{S^2_{\infty}({\Bbb H}^3)}
\newcommand{\PI}{\partial_{\infty}}
\newcommand{\SI}{S^2_{\infty}}
\newcommand{\BHH}{{\Bbb H}^3}
\newcommand{\CH}{\mathcal{C}(\Gamma)}
\newcommand{\BH}{\Bbb H}
\newcommand{\BR}{\Bbb R}
\newcommand{\BC}{\Bbb C}
\newcommand{\BZ}{\Bbb Z}

\begin{abstract}

Let $X$ be a Gromov hyperbolic $3$-space with cocompact metric, and $\Si$ be the sphere at infinity of $X$. We show that for any simple
closed curve $\Gamma$ in $\Si$, there exist a properly embedded least area plane in $X$ spanning $\Gamma$. This gives a positive answer to
Gabai's conjecture in [Ga1]. Soma has already proven this conjecture in [So1]. Our technique here is simpler and more general, and it can
be applied to many similar settings.

\end{abstract}

\section{Introduction}

Let $X$ be a Gromov hyperbolic $3$-space with cocompact metric, and $\Si$ be the sphere at infinity of $X$. Let $\Gamma$ be a given simple
closed curve in $\Si$. The asymptotic Plateau problem asks the existence of a least area plane $\Sigma$ in $X$ with asymptotic boundary
$\Gamma$, i.e. $\PI \Sigma = \Gamma$. The author gave a positive answer to this question, and showed the existence of such least area plane
asymptotic to given curve in $\Si$ in [Co].

Even though we know the existence of such least area planes in $X$ asymptotic to given curves in $\Si$, properly embeddedness of these
least area planes are also another interesting question came out in Gabai's renowned paper [Ga1]. In that paper, he solved the asymptotic
Plateau problem for hyperbolic $3$-spaces with cocompact metric, and he showed that any simple closed curve in $\SI(\widetilde{M})$, there
exists a least area plane spanning the given curve in $\widetilde{M}$ by using laminations where $M$ is a closed hyperbolic manifold with
any Riemannian metric. Also, he conjectured that there exist a $\emph{properly embedded}$ least area plane in $\widetilde{M}$ for a given
simple closed curve in $\SI(\widetilde{M})$, [Ga1, 3.12].

On the other hand, there is no known example of a nonproperly embedded least area plane in $\BH^3$ (or any Gromov hyperbolic $3$-space with
cocompact metric) spanning a simple closed curve in $\SI(\BH^3)$. There is a construction by Freedman and He communicated to Gabai [Ga2] of
a nonproper least area plane in $\BH^3$.  It is not clear to the author as to how a plane constructed in this manner can have limit set a
simple closed curve.

Teruhiko Soma gave a positive answer to the Gabai's conjecture in hyperbolic 3-space with cocompact metric in [So1]. Then, he generalized
his result to Gromov hyperbolic 3-spaces with cocompact metric in [So2]. In this paper, we will show the same results with a simpler and
more general technique which is applicable to many different settings. Our main result is the following:

\vspace{0.3cm}

\noindent {\bf Theorem 3.2.} Let $X$ be a Gromov hyperbolic $3$-space with cocompact metric, and $\Si$ be the sphere at infinity of $X$.
Let $\Gamma$ be a given simple closed curve in $\Si$. Then, there exist a properly embedded least area plane $\Sigma$ in $X$ with $\PI
\Sigma = \Gamma$.

\vspace{0.3cm}

A short outline of the proof is as follows. First, we consider the sequence of least area planes $\{\Sigma_i\}$ in $X$ with $\PI \Sigma_i =
\Gamma_i$ such that $\Gamma_i \to \Gamma$ in $\Si$ from one side. Then, we define a sequence of least area disks $\{D_i\}$ such that
$D_i\subset \Sigma_i$ and $\partial D_i \to \Gamma$. By using the results in [Co], we show that $\{D_i\}$ converges to a nonempty
lamination $\sigma$ by least area planes with $\PI \sigma = \Gamma$. Then, by using the fact that the least area disks stays in one side of
the lamination, we prove that this lamination contains only one least area plane $\mathcal{P}$, i.e. $\sigma = \mathcal{P}$. Since $\sigma$
is a closed subset of $X$, and so is the plane $\mathcal{P}$, then $\mathcal{P}$ must be a properly embedded least area plane with $\PI
\mathcal{P} = \Gamma$.

\subsection{Acknowledgements:}

I want to thank David Gabai and Yair Minsky for very helpful conversations.

\section{Preliminaries}

A $3$-manifold $M$ is called {\em Gromov hyperbolic manifold} if its fundamental group $\pi_1(M)$ is a word hyperbolic (or Gromov
hyperbolic) group [Gr]. We call $X$ as {\em Gromov hyperbolic $3$-space with cocompact metric} if $X$ is the universal cover of a
Riemannian closed orientable irreducible Gromov hyperbolic $3$-manifold $M$ where the metric on $X$ is induced by $M$. By [BM], $X$ is
homeomorphic to an open ball in $\BR^3$. Since $X$ is Gromov hyperbolic $3$-space, it has a natural compactification $\overline{M}$ where
$\overline{M}= M \cup \PI M$. Here, $\PI M$ is the sphere at infinity $\Si$, and a point on $\Si$ corresponds to an equivalence class of
infinite rays in $X$ where two rays are equivalent if they are asymptotic. A \textit{hyperbolic $3$-space with cocompact metric} is the
universal cover of a closed orientable irreducible hyperbolic $3$-manifold $M$ with any Riemannian metric $\rho$ on $M$. Hyperbolic
$3$-spaces with cocompact metric (in particular $\BH^3$) are special cases of Gromov hyperbolic $3$-spaces with cocompact metric.

We call a disk $D$ as a {\em least area disk} if $D$ has the smallest area among the disks with the same boundary $\partial D$. We call a
plane $P$ as a {\em least area plane} if any subdisk in the plane $P$ is a least area disk. A plane is \textit{properly embedded} if the
preimage of any compact set is compact in the plane.

A codimension-$1$ lamination $\sigma$ in $X$ is a foliation of a closed subset of $X$ with $2$-manifolds (\textit{leaves}) such that $X$ is
covered by charts of the form $I^2\times I$ where a leave passes through a chart in a slice of the form $I^2\times \{p\}$ for $p\in I$.
Here and later, we abuse notation by letting $\sigma$ also denote the underlying space of its lamination.

The sequence $\{D_{i}\}$ of embedded disks in a Riemannian manifold $X$ \textit{converges} to the lamination $\sigma$ if

i) For any convergent sequence $\{x_{n_i} \}$ in $X$ with  $x_{n_i} \in D_{n_i}$ where $n_i$ is a strictly increasing sequence, $\lim
x_{n_i} \in \sigma$.

ii) For any $x \in \sigma$, there exists a sequence $\{x_i\}$ with $x_i \in D_i$ and $\lim x_i = x$ such that there exist embeddings
$f_{i}: D^2 \to D_i$ which converge in the $C^{\infty }$-topology to a smooth embedding $f:D^2 \to L_{x}$, where $x_i \in f_{i}(Int(D^2))$,
and $L_{x}$ is the leaf of $\sigma $ through $x$, and $x\in f(Int(D^2))$.

We call such a lamination $\sigma$ as \textit{$D^2$-limit lamination}.

\begin{thm}\cite{Co}
Let $\Gamma $ be a simple closed curve in $\Si$ where $X$ is a Gromov hyperbolic $3$-space with cocompact metric.  Then there exists a
$D^2$-limit lamination $\sigma \subset X$ by least area planes spanning $\Gamma$.
\end{thm}

\begin{rmk}
Gabai proved this theorem for hyperbolic $3$-spaces with cocompact metric in [Ga1].
\end{rmk}

\section{Properly Embedded Least Area Plane}

In this section, we will prove the main theorem of the paper. First, we need a simple topological lemma.

\begin{lem} Let $\Sigma$ be a plane in $X$. If $\Sigma$ is closed in $X$, then $\Sigma$ is properly embedded.
\end{lem}

\begin{pf} Let $f: Int(D^2) \to X$ be an embedding with $f(Int(D^2))=\Sigma$. Let $K$ be a compact set in $X$. Since $\Sigma$ is closed, then $K \cap
\Sigma$ is also compact. Then, since $f$ is an embedding, $f^{-1}(K)=f^{-1}(K \cap \Sigma)$ is also compact. So, $\Sigma$ is a properly
embedded plane.
\end{pf}

Now, by the lemma, if a $D^2$-limit lamination $\sigma$ consists of only one plane $\Sigma$, i.e. $\sigma= \Sigma$, then $\Sigma$ must be
properly embedded as $\sigma$ is closed by definition. Hence, our aim is to construct a sequence of least area disks in $X$ converging to a
$D^2$-limit lamination $\sigma$ consisting of only one plane, which will be our properly embedded least area plane. Following lemma will be
very useful to construct such a sequence of least area disks in $X$.

\begin{lem}
Let $\Gamma_1$ and $\Gamma_2$ be two disjoint simple closed curves in $\Si$. If $\Sigma_1$ and $\Sigma_2$ are least area planes in $X$ with
$\PI \Sigma_i = \Gamma_i$, then $\Sigma_1$ and $\Sigma_2$ are disjoint, too.
\end{lem}

\begin{pf}
Assume that $\Sigma_1\cap\Sigma_2\neq\emptyset$. Since asymptotic boundaries $\Gamma_1$ and $\Gamma_2$ are disjoint, the intersection
cannot contain an infinite line. So, the intersection between $\Sigma_1$ and $\Sigma_2$ must contain a simple closed curve $\gamma$. Since
$\Sigma_1$ and $\Sigma_2$ are also minimal, the intersection must be transverse on a subarc of $\gamma$ by maximum principle.

Now, $\gamma$ bounds two least area disks $D_1$ and $D_2$ in $X$, with $D_i\subset\Sigma_i$. Now, take a larger subdisk $E_1$ of $\Sigma_1$
containing $D_1$, i.e. $D_1\subset E_1 \subset \Sigma_1$. By definition, $E_1$ is also an least area disk. Now, modify $E_1$ by swaping the
disks $D_1$ and $D_2$. Then, we get a new disk $E_1 '= \{E_1 - D_1\} \cup D_2$. Now, $E_1$ and $E_1 '$ have same area, but $E_1 '$ have a
folding curve along $\gamma$. By smoothing out this curve as in [MY], we get a disk with smaller area, which contradicts to $E_1$ being
least area.
\end{pf}

Now, we will prove the main theorem.

\begin{thm}Let $X$ be a Gromov hyperbolic $3$-space with cocompact metric, and $\Si$ be the sphere at infinity of $X$. Let $\Gamma$ be
a given simple closed curve in $\Si$. Then, there exist a properly embedded least area plane $\Sigma$ in $X$ with $\PI \Sigma = \Gamma$.
\end{thm}

\begin{pf}
Let $\Gamma$ be a simple closed curve in $\Si$. $\Gamma$ separates $\Si$ into two parts, say $\Omega^+$ and $\Omega^-$ which are open
disks. Define sequences of pairwise disjoint simple closed curves $\{\Gamma_i\}$  such that $\Gamma_i = \partial E_i$ where $E_i$ is a
closed disk in $\Omega^+$ for any $i$, and $E_i\subset Int(E_j)$ for any $i<j$. Moreover, $\Omega^+ =\bigcup_i E_i$ and $\Gamma_i
\rightarrow \Gamma$.

By Theorem 2.1, for any $\Gamma_i \subset \Si$, there exist a least area plane $\Sigma_i$ with $\PI \Sigma_i = \Gamma_i$. Note that these
least area planes are separating in $X$ by their construction in Theorem 2.1. Also, let $p^+ \in Int(E_1) \subset \Omega^+$ and $p^- \in
\Omega^-$ be two points belonging different components of $\Si-\Gamma$. Let $\beta$ be an infinite geodesic in $X$ asymptotic to $p^+$ and
$p^-$. Since each $\overline{\Sigma_i}$ is separating in $\overline{X}$, $\beta$ intersects $\Sigma_i$ for any $i$. Let $x_i$ be a point in
$\Sigma_i \cap \beta$ for any $i$.

Now, we define the sequence of least area disks. Let $D_i$ be a least area disk in $\Sigma_i$ containing $x_i$ such that $i<
d_X(x_i,\partial D_i) < i+1$ where $d_X$ is the extrinsic distance in $X$. Alternatively, by using a result in [An], one can define $D_i$
as follows. $\Sigma_i \cap B_r(x_i)$ is a collection of disjoint disks for almost all $r$, where $B_r(x_i)$ is the $r$-ball in $X$ with
center $x_i$. Then, define $D_i$ is the component containing $x_i$ in $B_r(x_i)\cap \Sigma_i$ where $r_i$ is a regular value between $i$
and $i+1$.

We claim that the sequence of least area disks $\{D_i\}$ in $X$ converge to a nonempty lamination $\sigma$.  As $\alpha_i =\partial D_i$
converges to $\Gamma$, by Theorem 2.1, $D_i$ converges to a lamination $\sigma$ (possibly empty) by least area planes with $\PI \sigma
=\Gamma$. Now, we show that $\sigma$ is a nonempty lamination by least area planes. Since $p^+ \in \Omega^+$ and $p^-\in \Omega^-$ where
$\Omega^\pm$ are open disks in $\Si$, we can find a sufficiently small $\epsilon$ such that $B_\epsilon(p^+) \subset Int(E_1)$ and
$B_\epsilon(p^-) \subset \Omega^-$. Let $\gamma^+ = \partial B_\epsilon(p^+)$ and $\gamma^- = \partial B_\epsilon(p^-)$. Then, $\gamma^\pm$
are two simple closed curves disjoint from $\Gamma$. By Theorem 2.1, there are least area planes $P^+$ and $P^-$ in $X$ with $\PI P^\pm =
\gamma^\pm$. By Lemma 3.2, $P^\pm$ are disjoint from $\Sigma_i$ for each $i$. So, the finite segment $\hat{\beta}$ of $\beta$ between $P^+$
and $P^-$ contains all intersection points $\Sigma_i \cap \beta$ for any $i$. Hence, $\{x_i\}\subset\hat{\beta}$. Since $\hat{\beta}$ is
finite segment, it is compact, and the sequence $\{x_i\}$ has a limit point $x$. This shows that the sequence $\{D_i\}$ has nonempty limit.
So, $\sigma$ is a nonempty lamination by least area planes with $\PI \sigma = \Gamma$ as claimed.

Now, we want to show that $\sigma$ consists of only one least area plane $\mathcal{P}$, i.e. $\sigma = \mathcal{P}$. Assume that $\sigma$
contains more than one least area plane. Recall that $\sigma$ is a collection of disjoint least area planes asymptotic to $\Gamma$, and
$\sigma$ is a closed subset of $X$. Note that each plane $\mathcal{P}$ in $\sigma$ is separating in $X$ by construction. Consider the
components of $X- \sigma$. Let $X^+$ be the component of $X - \sigma$ with $\PI X^+ = \Omega^+$. Note that as $\sigma$ is closed, each
component of $X-\sigma$ is open, and so is $X^+$. Since for each $i$, $D_i \subset \Sigma_i$ and $\PI \Sigma_i = \Gamma_i \subset
\Omega^+$, by Lemma 3.2, $D_i \subset X^+$. Since each plane $\mathcal{P}$ in $\sigma$ is separating, there is a plane $\mathcal{P}_1$ such
that $\mathcal{P}_1 = \partial X^+$. Let $\mathcal{P}_2$ be another least area plane in $\sigma$. We claim that there cannot be such a
plane $\mathcal{P}_2$ because of the special properties of the sequence of least area disks $\{D_i\}$. Let $Y^+$ be the component of
$X-\mathcal{P}_2$ with $\PI Y^+ = \Omega^+$. Clearly, $X^+ \subset Y^+$. Moreover, since either there exist an open complementary region (a
component of $X-\sigma$), or there is an open region foliated by least area planes in $\sigma$ between $\mathcal{P}_1$ and $\mathcal{P}_2$,
hence $X^+ \subsetneq Int(Y^+)$. This means $\mathcal{P}_1$ forms a barrier between the sequence of least area disks $\{D_i\}$ and
$\mathcal{P}_2$. In other words, the sequence $\{D_i\}$ cannot reach $\mathcal{P}_2$, so $\mathcal{P}_2$ cannot be in the limit. More
precisely, if $q \in \mathcal{P}_2$, since $D_i\subset X^+$ for any $i$, and $X^+ \subsetneq Int(Y^+)$, there is no sequence $\{q_i\}$ with
$q_i \in D_i$ such that $q_i \to q$. So, $\mathcal{P}_2$ cannot be in the limit of $\{D_i\}$. Hence, this shows $\sigma$ consists of only
one plane $\mathcal{P}_1$, i.e. $\sigma = \mathcal{P}_1$.

As $\sigma= \mathcal{P}_1$, and $\sigma$ is a lamination, $\mathcal{P}_1$ is closed. By Lemma 3.1, $\mathcal{P}_1$ is properly embedded in
$X$. This shows $\mathcal{P}_1$ is a properly embedded least area plane in $X$ with $\PI \mathcal{P}_1 = \Gamma$ as desired.
\end{pf}

\section{Final remarks}

In this paper, we showed that for any simple closed curve in asymptotic boundary of a Gromov hyperbolic $3$-space with cocompact metric,
there exist a properly embedded least area plane in the space spanning the curve. The technique we use here is very general, and can be
applied to many similar settings. For example, by using results of Lang in [La], these results can be extended to Gromov hyperbolic
Hadamard $3$-spaces. Similarly, if we have a positive solution to the asymptotic Plateau problem in a $3$-space $X$, then the methods in
this paper can naturally be generalized to this case.

On the other hand, as we point out in the introduction, while we prove the existence of properly embedded least area planes in hyperbolic
$3$-spaces, there is no known nonproperly embedded least area plane in $\BH^3$. By using the regularity results of Hardt and Lin in [HL],
we know that if a simple closed curve in $\SI(\BH^3)$ is $C^1$, then any least area plane asymptotic to the curve is properly embedded.
However, this is not known yet for $C^0$ simple closed curves in $\SI(\BH^3)$. Intuitively, being least area and nonproperly embeddednes
are contradicting notions. As, nonproper embeddednes must produce monogons, one can get a contradiction by using Hass and Scott's surgery
arguments for least area objects in [HS]. Hence, we have the following conjecture.

\vspace{0.3cm}

\noindent {\bf Conjecture:} Let $\Gamma$ be a simple closed curve in $\SI (\BH^3)$. Then, any least area plane $\Sigma$ in $\BH^3$ with
$\PI \Sigma = \Gamma$ is properly embedded.

\end{document}